%% file: main.tex
\documentclass[smallextended,envcountsect]{svjour3}
\smartqed 
\usepackage{graphicx}

\usepackage{amsmath,amsfonts,amssymb,mathtools}
\usepackage{stmaryrd}
\usepackage{algorithm}
\usepackage{algpseudocode}
\usepackage{siunitx}
\usepackage{xfrac}
\usepackage{xcolor}
\usepackage{hyperref}
\usepackage{soul}
\usepackage{booktabs,multirow,multicol,arydshln}
\usepackage[noabbrev,capitalise]{cleveref}
\crefname{equation}{}{}
\colorlet{inlinkcolor}{green!50!black}
\colorlet{exlinkcolor}{red!50!black}
\hypersetup{
colorlinks = true,
allcolors = inlinkcolor,
urlcolor = exlinkcolor,
}
\input{defs}

\begin{document}

\title{Computing Conjugate Barrier Information for Nonsymmetric Cones}


\author{Lea Kapelevich   \and  Erling D. Andersen \and Juan Pablo Vielma}

\institute{Lea Kapelevich, Corresponding author \at
        Operations Research Center, MIT \\
        Cambridge, MA, USA \\
        lkap@mit.edu   
        \and
        Erling D. Andersen \at
        MOSEK ApS \\
        Copenhagen, Denmark \\
        e.d.andersen@mosek.com
        \and
        Juan Pablo Vielma  \at
        Google Research and MIT Sloan School of Management \\
        Cambridge, MA, USA \\
        jvielma@google.com
}

\date{Received: date / Accepted: date}

\maketitle

\begin{abstract}
\input{abstract}
\end{abstract}
\\
\\
{
This preprint has not undergone all peer review and post-submission improvements or corrections. 
The Version of Record of this article is published in the Journal of Optimization Theory and Applications, and is available online at \url{https://doi.org/10.1007/s10957-022-02076-1}.
\\
}

\input{introduction}

\input{preliminaries}
\section{Conjugate Gradients}

In \crefrange{sec:sumlog}{sec:inf} we offer efficient procedures for evaluating conjugate gradients.
We defer some derivations to \cref{sec:conjproofs} to ease readability.
\subsection{Logarithm Cone and Log-determinant Cone}
\label{sec:sumlog}

Let $\varphi: \bbR^d_{++} \to \bbR$ be the function $\varphi(w) \coloneqq \tsum_{i \in \iin{d}} \log w_i$.
Define the \emph{logarithm} cone:
\begin{equation*}
    \klog \coloneqq \cl \bigl\{ (u, v, w) \in \bbR \times \bbR_{++} \times \bbR_{++}^d :  u \leq v \varphi \bigl( \tfrac{w}{v} \bigr) \bigr\},
\end{equation*}
which admits the $(2 + d)$-LHSCB \cite[Section 6]{coey2021conic}:
\begin{equation*}
    f(u, v, w) =  -\log \bigl(v \varphi \bigl( \tfrac{w}{v} \bigr) - u \bigr) - \log(v) - \tsum_{i \in \iin{d}} \log(w_i).
\end{equation*}
The dual cone is given by \cite[Section 5.5]{coey2020solving}:
\begin{equation*}
    \klog^\ast \coloneqq \cl \bigl\{ (p, q, r) \in \bbR_{--} \times \bbR \times \bbR^d_{++} : q \geq p \tsum_{i \in \iin{d}} \log \bigl(-\tfrac{r_i}{p} \bigr) + p d \bigr\}.
\end{equation*}
Let $\omega$ denote the Wright omega function \cite{corless2002wright}, which can be well approximated in $O(1)$ time and satisfies for all real $\beta$:
\begin{equation}
    \omega(\beta) + \log(\omega(\beta)) = \beta.
    \label{eq:omegachar}
\end{equation}
The Wright omega function is used in \cite[Chapter 8]{serrano2015algorithms} for deriving the conjugate barrier of a three-dimensional variant of $\klog$ and in \cite[Section 7.1]{nemirovski2005cone} for the conjugate barrier of the epigraph of the univariate entropy function.
\begin{proposition}
\label{prop:sumlog}
The conjugate gradient at $(p, q, r) \in \intr(\klog^\ast)$ has components:
\begin{subequations}
\label{eq:tgsumlog}
\begin{align}
    \tg_p &= \tfrac{-d - 2 + \sfrac{q}{p} + 2 \bomega}{p (1 - \bomega)},
    \\
    \tg_q &= -\tfrac{1}{p (1 - \bomega)},
    \\
    \tg_{r_i} &= \tfrac{\bomega}{r_i(1 - \bomega)} & \forall i \in \iin{d},
    \label{eq:tgsumlog:ri}
\end{align}
\end{subequations}
where, $\bar{\omega} \coloneqq d \cdot \omega \bigl( \sfrac{1}{d} \bigl( 1 + d - \sfrac{q}{p} + \tsum_{i \in \iin{d}} \log(-\sfrac{r_i}{p}) \bigr) - \log(d) \bigr)$.
\end{proposition}
The proof is given in \cref{sec:conjproofs:sumlog}.
By substituting \cref{eq:tgsumlog} in \cref{eq:fast2}, we obtain the conjugate barrier:
\begin{equation*}
    f^\ast(p, q, r) = 
    -2 - d - 2 \log(-p) - \log \Bigl( \tfrac{(\bomega - 1)^{d+1}}{\bomega^d} \Bigr) - \tsum_{i \in \iin{d}} \log(r_i).
\end{equation*}
The conjugate gradient for $\klog$ can be easily modified to obtain a conjugate gradient for the \emph{log-determinant cone} \cite[Section 5.6]{coey2020solving}:
\begin{equation*}
    \klogdet \coloneqq \cl \bigl\{ (u, v, W) \in \bbR \times \bbR_{++} \times \bbS^d_{++} :  u \leq v \logdet (\sfrac{W}{v} ) \bigr\},
\end{equation*}
which has the dual:
\begin{equation*}
    \klogdet^\ast \coloneqq \cl \bigl\{ (p, q, R) \in \bbR_{--} \times \bbR \times \bbS^d_{++} :  q \geq p \logdet (-\sfrac{R}{p} ) + p d \bigr\}.
\end{equation*}
Let $W = U_W \Diag(\lambda_W) U_W^\top$ be the eigendecomposition of $W$.
$\klogdet$ admits the $(2+d)$-LHSCB \cite[Section 6]{coey2021conic}:
\begin{equation*}
    F(u, v, W) = -\log ( v \varphi (\sfrac{\lambda_W}{v} ) - u ) - \log(v) - \tsum_{i \in \iin{d}} \log(\lambda_{W,i}) = f(u, v, \lambda_W).
\end{equation*}
\begin{proposition}
\label{prop:logdet}
Let $(p, q, R) \in \intr(\klogdet^\ast)$, let $R = U_R \Diag(\lambda_R) U_R^\top$ be the eigendecomposition of $R$, and let $\bOmega \coloneqq d \cdot \omega \bigl( \sfrac{1}{d} \bigl( 1 + d - \sfrac{q}{p} + \tsum_{i \in \iin{d}} \log(-\sfrac{\lambda_{R,i}}{p}) \bigr) - \log(d) \bigr)$.
The conjugate gradient $G^\ast$ has components:
\begin{align*}
    G^\ast_p &= \tg_p(p, q, \lambda_R) = \tfrac{-d - 2 + \sfrac{q}{p} + 2 \bOmega}{p (1 - \bOmega)},
    \\
    G^\ast_q &= \tg_q(p, q, \lambda_R) = -p^{-1} (1 - \bOmega)^{-1},
    \\
    G^\ast_R &= U_R \Diag(\tg_r(p, q, \lambda_R)) U_R^\top = \tfrac{\bOmega}{(1 - \bOmega)} R^{-1}.
\end{align*}
\end{proposition}
{\it Proof}
For fixed $u$ and $v$, $F$ is a unitarily invariant function of $W$.
Due to \cref{eq:unitarily:G}, the gradient of $F$ is $G_u = g_u(u, v, \lambda_W)$, $G_v = g_v(u, v, \lambda_W)$, and $G_W = U_W \Diag(g_w(u, v, \lambda_W)) U_W^\top$.
The result can be verified from this and \cref{eq:grbilinear}.
\qed
\subsection{Hypograph Power Cone and Root-determinant Cone}
\label{sec:hpower}

Let $\varphi: \bbR^d \to \bbR$ be the function%
\footnote{%
We reuse symbols with similar roles across subsections; their meaning should be taken from the definition within each subsection.
}
$\varphi(w) \coloneqq \tprod_{i \in \iin{d}} w_i^{\alpha_i}$, parametrized by $\alpha = (\alpha_1, \ldots, \alpha_d)$ such that $\langle e, \alpha \rangle = 1$ and $\alpha > 0$.
Define the \emph{hypograph-power} cone:
\begin{equation*}
    \khpower \coloneqq \bigl\{ (u, w) \in \bbR \times \bbR_+^d : u \leq \varphi(w) \bigr\},
\end{equation*}
which admits the $(1+d)$-LHSCB \cite[Section 5.4.7]{nesterov2018lectures}:
\begin{equation}
    f(u, w) = -\log(\varphi(w) - u) - \tsum_{i \in \iin{d}} \log(w_i).
    \label{eq:hpower:f}
\end{equation}
The cone $\khpower$ is the \emph{power mean cone} in \cite{coey2021performance}.
In the special case where $\alpha = \sfrac{e}{d}$, we call $\khpower$ the \emph{hypograph geometric mean} cone, $\khgeom$.
The dual cone is \cite{coey2021hypatia}:
\begin{equation*}
   \khpower^\ast \coloneqq \bigl\{ (p, r) \in \bbR_{-} \times \bbR_+^d : -p \leq \varphi \bigl(\tfrac{r}{\alpha} \bigr) \bigr\}.
   \label{eq:hpowerdual}
\end{equation*}
\begin{lemma}
\label{lem:hpower}
Let $(p, r) \in \intr(\khpower^\ast)$ parametrized by $\alpha$.
The unique root of $h(y) \coloneqq \tsum_{i \in \iin{d}} \alpha_i \log (y - p \alpha_i) - \log(\varphi(r))$ can be approximated by a quadratically convergent Newton-Raphson method starting from $0$.
\end{lemma}
The proof is given in \cref{sec:conjproofs:hpower}.
\begin{proposition}
\label{prop:hpower}
The conjugate gradient at $(p, r) \in \intr(\khpower^\ast)$ has components:
\begin{subequations}
\label{eq:khpower:tg}
\begin{align}
    \tg_p &= -p^{-1} - {\hat{y}}^{-1},
    \\
    \tg_{r_i} &= \tfrac{p \alpha_i {\hat{y}}^{-1} - 1}{r_i} & \forall i \in \iin{d},
\end{align}
\end{subequations}
where $\hat{y}$ is the root of $h$ from \cref{lem:hpower}.
In the case where $\alpha = \sfrac{e}{d}$, the conjugate gradient at $(p, r) \in \intr(\khgeom^\ast)$ can be written more simply:
\begin{subequations}
\begin{align}
    \tg_p &= -p^{-1} - (\varphi(r) + \sfrac{p}{d})^{-1},
    \\
    \tg_{r_i} &= -\tfrac{\varphi(r)}{r_i (\varphi(r) + \sfrac{p}{d})} & \forall i \in \iin{d}.
\end{align}
\label{eq:khgeom:tg}
\end{subequations}
\end{proposition}
The proof is given in \cref{sec:conjproofs:hpower}.
Substituting \cref{eq:khgeom:tg} in \cref{eq:fast2}, we obtain a simple expression for the conjugate barrier of $\khgeom$:
\begin{equation*}
    f^\ast(p, r) = -1 - d - d \log \big( \tfrac{d  \varphi(r) + p}{d  \varphi(r)} \big) -
    \log (-p) - \tsum_{i \in \iin{d}} \log \big( r_i \big).
    \label{eq:hpower:fast}
\end{equation*}
The conjugate gradient for $\khgeom$ can be easily modified to obtain a conjugate gradient for the \emph{root-determinant cone} \cite[Section 5.4]{coey2020solving}:
\begin{equation*}
    \krtdet \coloneqq \bigl\{ (u, W) \in \bbR \times \bbS^d_+ : u \leq \det(W)^{1/d} \bigr\},
\end{equation*}
which has the dual:
\begin{equation*}
    \krtdet^\ast \coloneqq \bigl\{ (p, R) \in \bbR_{-} \times \bbS^d_+ : -p \leq d \det(R)^{1/d} \bigr\}.
\end{equation*}
Let $W = U_W \Diag(\lambda_W) U_W^\top$ be the eigendecomposition of $W$.
$\krtdet$ admits the $(1+d)$-LHSCB \cite{coey2021conic}:
\begin{equation*}
    F(u, W) = -\log(\varphi(\lambda_W) - u) - \tsum_{i \in \iin{d}} \log(\lambda_{W,i}) = f(u, \lambda_W).
\end{equation*}
\begin{proposition}
Let $(p, R) \in \intr(\krtdet^\ast)$ and $R = U_R \Diag(\lambda_R) U_R^\top$ be the eigendecomposition of $R$.
The conjugate gradient $G^\ast$ has components:
\begin{align*}
    G^\ast_p &= g_p^\ast(p, \lambda_R) = -p^{-1} - (\det(R)^{1/d} + \sfrac{p}{d})^{-1},
    \\
    G^\ast_R &= U_R \Diag(g_r^\ast(p, \lambda_R)) U_R^\top = -\tfrac{\det(R)^{1/d}}{\det(R)^{1/d} + \sfrac{p}{d}} R^{-1}.
\end{align*}
\end{proposition}
{\it Proof}
For fixed $u$, $F$ is a unitarily invariant function of $W$.
Similar to \cref{prop:logdet}, the result can be verified from \cref{eq:unitarily:G,eq:grbilinear}.
\qed
\subsection{Radial Power Cone}
\label{sec:rpower}

Let $\varphi: \bbR^d \to \bbR$ be the function $\varphi(w) \coloneqq \tprod_{i \in \iin{d}} w_i^{2 \alpha_i}$, parametrized by $\alpha = (\alpha_1, \ldots, \alpha_d)$ such that $\langle e, \alpha \rangle = 1$ and $\alpha > 0$.
Define the \emph{radial-power} cone:
\begin{equation*}
    \krpower \coloneqq \bigl\{ (u, w) \in \bbR^{d_1} \times \bbR_+^{d_2} : \lVert u \rVert \leq \tprod_{i \in \iin{d_2}} w_i^{\alpha_i} = \sqrt{\varphi(w)} \bigr\},
\end{equation*}
which admits the $(1+d_2)$-LHSCB \cite[Theorem 1]{roy2021self}:
\begin{equation*}
    f(u, w) = -\log(\varphi(w) - \lVert u \rVert^2) - \tsum_{i \in \iin{d_2}} (1 - \alpha_i) \log(w_i).
\end{equation*}
This is the \emph{generalized power cone} in \cite{coey2021performance,chares2009cones,roy2021self}.
Note that $f$ is not equivalent to the barrier from \cref{eq:hpower:f}, even when $d_1 = 1$.
Hence the conjugate barrier and its derivatives take different forms from our results in \cref{sec:hpower}.
In the special case where $\alpha = \sfrac{e}{d_2}$ and $d_1 = 1$, we call $\krpower$ the \emph{radial geometric mean} cone $\krgeom$.
The dual cone is given by \cite[Theorem 4.3.1]{chares2009cones}:
\begin{equation*}
   \krpower^\ast \coloneqq \bigl\{ (p,r) \in \bbR^{d_1} \times \bbR_+^{d_2} : \lVert p \rVert \leq \tprod_{i \in \iin{d_2}} \bigl(\tfrac{r_i}{\alpha_i} \bigr)^{\alpha_i} \bigr\}.
\end{equation*}

\begin{lemma}
\label{lem:rpower}
Let $(p, r) \in \intr(\khpower^\ast)$ parametrized by $\alpha$ and $p > 0$.
The unique positive root of $h(y) \coloneqq \tsum_{i \in \iin{d_2}} 2 \alpha_i \log (2 \alpha_i y^2 + \sfrac{2 y (1 + \alpha_i)}{p} ) - \log(\varphi(r)) - \log ( \sfrac{2 y}{p} + y^2 ) - 2 \log ( \sfrac{2 y}{p} ) $
can be approximated by a quadratically convergent Newton-Raphson method starting from $ y_{-} \coloneqq -p^{-1} + d_2 \tfrac{p + \sqrt{\varphi(r) (\sfrac{d_2^2}{p^2} \varphi(r) + d_2^2 - 1)}}{\varphi(r) d_2^2 - p^2}$.
\end{lemma}
The proof is given in \cref{sec:conjproofs:rpower}.
\begin{proposition}
\label{prop:rpower}
The conjugate gradient at $(p, r) \in \intr(\krpower^\ast)$ is given by:
\begin{subequations}
\label{eq:rpower:tg}
\begin{align}
    \tg_{p_i} &= 
    1_{p \neq 0} \cdot
    \hat{y} \tfrac{p_i}{\lVert p \rVert},
    & \forall i \in \iin{d_1},
    \\
    \tg_{r_i} &=  -\tfrac{\alpha_i (1 + p \hat{y}) + 1}{r_i} & \forall i \in \iin{d_2},
    \label{eq:rpower:tgri}
\end{align}
\end{subequations}
where $\hat{y}$ is the positive root of $h$ from \cref{lem:rpower}.
In the case where $d_1 = 1$ and $\alpha = \sfrac{e}{d_2}$, the conjugate gradient at $(p, r) \in \intr(\krgeom^\ast)$ is:
\begin{subequations}
\begin{align}
    \tg_{p} &= y_{-},
    \\
    \tg_{r_i} &= -r_i^{-1} \Bigl( \tfrac{p^2 + \sqrt{\varphi(r) (d_2^2 \varphi(r) + d_2^2 p^2 - p^2)}}{\varphi(r) d_2^2 - p^2} + 1 \Bigr) & \forall i \in \iin{d_2}.
    \label{eq:krgeom:tgr}
\end{align}
\label{eq:krgeom:tg}
\end{subequations}
\end{proposition}
Note $y_{-}$ is from \cref{lem:rpower}.
The proof is given in \cref{sec:conjproofs:rpower}.
\subsection{Infinity Norm Cone and Spectral Norm Cone}
\label{sec:inf}

Define the \emph{infinity norm cone} \cite[Section 5.1]{coey2020solving}:
\begin{equation*}
    \kinf \coloneqq \bigl\{ (u, w) \in \bbR \times \bbR^d : u \geq \Vert w \Vert_{\infty} \bigr\},
\end{equation*}
which admits the $(1+d)$-LHSCB \cite[section 7.5]{guler1996barrier}:
\begin{equation*}
    f(u, w) = -\tsum_{i \in \iin{d}} \log(u^2 - w_i^2) + (d - 1) \log(u).
\end{equation*}
The dual cone is the epigraph of the $\ell_1$ norm function:
\begin{equation*}
    \kinf^\ast \coloneqq \bigl\{ (p, r) \in \bbR \times \bbR^d : p \geq \Vert r \Vert_{1} \bigr\}.
\end{equation*}
\begin{lemma}
\label{lem:inf}
Let $(p, r) \in \intr(\kinf^\ast)$.
The unique negative root of $ h(y) \coloneqq p y + \tsum_{i \in \iin{d}} \sqrt{1 + r_i^2 y^2} + 1$ can be approximated by a quadratically convergent Newton-Raphson method starting from $\max \{ -(p - \lVert r \rVert_1)^{-1}, -\sfrac{(d + 1)}{p}\}$.
\end{lemma}
The proof is given in \cref{sec:conjproofs:hpower:inf}.
\begin{proposition}
\label{prop:inf}
The conjugate gradient at $(p, r) \in \intr(\kinf^\ast)$ is given by:
\begin{subequations}
\begin{align}
    \tg_p &= \hat{y},
    \\
    \tg_{r_i} &= 
    1_{r_i \neq 0} \cdot \tfrac{\sqrt{1 + {\hat{y}}^2 r_i^2} - 1}{r_i},
    & \forall i \in \iin{d},
    \label{eq:inf:tgr}
\end{align}
\end{subequations}
where $\hat{y}$ is the negative root of $h$ from \cref{lem:inf}.
\end{proposition}
The proof is given in \cref{sec:conjproofs:hpower:inf}.
The conjugate gradient for $\kinf$ can be easily modified to obtain a conjugate gradient for the \emph{spectral norm cone} \cite[Section 5.2]{coey2020solving}:
\begin{equation*}
    \kspec \coloneqq \{ (u, W) \in \bbR \times \bbR^{d_1 \times d_2} : u \geq \sigma_{\max}(W) \},
\end{equation*}
where $\sigma_{\max}$ is the maximum singular value function and $d_1 \leq d_2$.
Let $R = U_R \Diag(\sigma_R) V_R^\top$ be the singular value decomposition of $R \in \bbR^{d_1 \times d_2}$.
The dual cone is the epigraph of the nuclear norm:
\begin{equation*}
    \kspec^\ast \coloneqq \{ (p, R) \in \bbR \times \bbR^{d_1 \times d_2} : p \geq \tsum_{i \in \iin{d_1}} \sigma_{R,i} \}.
\end{equation*}
Let $W = U_W \Diag(\sigma_W) V_W^\top$ be the singular value decomposition of $W$.
$\kspec$ admits the $(1+d_1)$-LHSCB:
\begin{equation*}
    F(u, W) = -\tsum_{i \in \iin{d_1}} \log(u^2 - \sigma_{W,i}^2) + (d_1 - 1) \log(u) = f(u, \sigma_W).
    \label{eq:inf:F}
\end{equation*}
\begin{proposition}
The conjugate gradient $G^\ast$ at $(p, R) \in \intr(\kspec^\ast)$ has components:
\begin{align*}
    G^\ast_p &= \tg_p(p, \sigma_R),
    \\
    G^\ast_R &= U_R \Diag(\tg_r(p, \sigma_R)) V_R^\top.
\end{align*}
\end{proposition}
{\it Proof}
For fixed $u$, $F$ is a unitarily invariant function of $W$.
Similar to \cref{prop:logdet}, the result can be verified using \cref{eq:unitarily:G,eq:grbilinear}.
\qed
\section{Inverse Hessians}
\label{sec:invhess}

In \cite{coey2021conic} the authors derive efficient inverse Hessian operators for a number of cones, including $\klog$, $\klogdet$, $\khgeom$, and $\krtdet$.
Inverse Hessians for $\kinf$ and $\kspec$ are described in \cite{coey2022thesis}.
The authors are motivated by the need to apply inverse Hessians in several parts of the algorithm implemented by the Hypatia solver \cite{coey2021performance}.
The inverse Hessian is used to measure proximity to the central path as well as some of the linear system solving methods described in \cite[Appendix A]{coey2021performance}.
Our derivations of the conjugate gradients offer an alternative method for deriving inverse Hessian operators.
Since these have not been previously written for $\khpower$ or $\krpower$, we derive those inverse Hessian operators here.

A general technique for obtaining an inverse Hessian operator using an oracle for the conjugate gradient may be described as follows.
For $\K, \K^\ast \subset \bbR^{1 + d}$, where $\K$ is either $\khpower$ or $\krpower$ parametrized by $\alpha$, let $\tu = (u, w) \in \intr(\K)$, and $\tx = (x, z) \in \bbR^{1 + d}$ be an arbitrary direction.
Due to \cite[Equation (2.11)]{nesterov1997self}:
\begin{equation}
    H(\tu)^{-1} \cdot \tx = H^\ast(-g(\tu)) \cdot \tx.
    \label{eq:dirHfact}
\end{equation}
Therefore to derive the inverse Hessian operator $H(\tu)^{-1} \cdot \tx$, we may derive an expression for $H^\ast(\tp)\cdot \tx$, for arbitrary $\tp = (p, r) \in \intr(\K^\ast)$, and substitute $-g(\tu)$ for $\tp$.
In a practical implementation of a PDIPM, $g(\tu)$ is usually already available at the time when $H(\tu)^{-1}$ is needed.
Note that:
\begin{equation*}
    H^\ast(\tp) \cdot \tx = \tfrac{\diff}{\diff t} g^\ast(\tp + t \tx) \big\rvert_{t=0}.
\end{equation*}
The right hand side can be deduced from our calculations for $\tg$.
\subsection{Hypograph Power Cone}
\label{sec:invhess:hpower}

\begin{proposition}
\label{prop:invhess:hpower}
The inverse Hessian operator at $\tu$ in the direction $\tx$ is:
\begin{subequations}
\label{eq:hinvhpower}
\begin{align}
    (H(\tu)^{-1}\cdot \tx)_u &= \bigl( (\varphi(w) - u)^2 + \tfrac{k_2}{k_3} u^2 \bigr) x - \tfrac{\varphi(w)}{k_3} \langle z, \tfrac{\alpha w}{k_1}  \rangle,
    \\
    (H(\tu)^{-1}\cdot \tx)_{w_i} &= \tfrac{w_i^2}{k_{1,i}} z_i + \tfrac{\alpha_i w_i}{k_{1,i}} \tfrac{\varphi(w)}{k_3} x + \tfrac{g_u \varphi(w)}{k_3} \langle z, \tfrac{\alpha w}{k_1}  \rangle \tfrac{\alpha_i w_i}{k_{1,i}} & \forall i \in \iin{d},
\end{align}
\end{subequations}
where $\varphi$ is defined as in \cref{sec:hpower}, and:
\begin{equation*}
    k_{1,i} \coloneqq 1 + \alpha_i \varphi(w) g_u \quad \forall i \in \iin{d}, \quad
    k_2 \coloneqq \tsum_{i \in \iin{d}} \tfrac{\alpha_i^2}{k_{1,i}}, \quad
    k_3 \coloneqq 1 - \varphi(w) g_u k_2.
\end{equation*}
\end{proposition}
The proof is given in \cref{sec:hessproofs:hpower}.
\subsection{Radial Power Cone}
\label{sec:invhess:rpower}

\begin{proposition}
\label{prop:invhess:rpower}
The inverse Hessian operator at $\tu$ in the direction $\tx$ is:
\begin{subequations}
\label{eq:Hinvrpower}
\begin{align}
    (H(\tu)^{-1} \cdot \tx)_{u_i} &= \tfrac{\zeta x_i}{2} + \tfrac{u_i}{k_3} \big( \tfrac{\zeta k_3 - 2 k_2 \varphi(w)}{k_1} \langle x, u \rangle - \langle \sfrac{\alpha}{g_w}, z \rangle \big) & \forall i \in \iin{d_1},
    \\
    (H(\tu)^{-1} \cdot \tx)_{w_i} &= -\tfrac{w_i z_i}{g_{w_i}} - \tfrac{\alpha_i}{k_3 g_{w_i}} \big(\langle x, u \rangle - \tfrac{2 \lVert u \rVert^2}{\zeta} \langle \sfrac{\alpha}{g_w}, z \rangle \big) & \forall i \in \iin{d_2},
\end{align}
\end{subequations}
where $\varphi$ defined as in \cref{sec:rpower}, $\zeta \coloneqq \varphi(w) - \lVert u \rVert^2$, and:
\begin{equation*}
    k_1 \coloneqq \varphi(w) + \lVert u \rVert^2, \quad
    k_2 \coloneqq \langle \sfrac{\alpha}{w}, \sfrac{\alpha}{g_w} \rangle , \quad
    k_3 \coloneqq \tfrac{k_1}{2 \varphi(w)} + 2 k_2 \tfrac{\lVert u \rVert^2}{\zeta}.
\end{equation*}
\end{proposition}
The proof is given in \cref{sec:hessproofs:rpower}.
\section{Proofs of Conjugate Gradients}
\label{sec:conjproofs}

For convenience, we adapt \cite[Theorem 1.9]{suli2003introduction} as a Lemma below (without proof), which offers some conditions for quadratic convergence of the Newton-Raphson method.
\begin{lemma}
\label{lem:lr}
Suppose that the continuous, real-valued function $h$ has a continuous second derivative on the closed interval $I_\delta \coloneqq [\hat{y} - \delta, \hat{y} + \delta]$, $\delta > 0$, such that $h(\hat{y}) = 0$ and $h''(\hat{y}) \neq 0$.
Suppose further that $\lvert h''(x) \rvert / \lvert h'(y) \rvert$ is bounded for all $x, y \in I_\delta$.
If there exists a real number $y_+ > \hat{y}$ (or $y_{-} < \hat{y}$), such that in the interval $J \coloneqq [\hat{y}, y_+]$ (or $J \coloneqq [y_{-}, \hat{y}]$) both $h'$ and $h''$ are positive (or one is negative), then the Newton-Raphson method converges quadratically to $\hat{y}$ from any starting value in $J$.
\end{lemma}
\cref{lem:lr} is used in our proofs of \cref{lem:hpower,lem:rpower,lem:inf}.
\subsection{Logarithm Cone}
\label{sec:conjproofs:sumlog}

{\it Proof of \cref{prop:sumlog}}
For convenience, let $\zeta$ be the function $\zeta(u, v, w) \coloneqq v \varphi(\sfrac{w}{v}) - u$, where $\varphi$ is from \cref{sec:sumlog}.
Let $(u, v, w) \in \intr(\klog)$.
Then the gradient of $f$ with respect to components, $u$, $v$, and $w$ is:
\begin{subequations}
\begin{align}
g_u &= \zeta(u, v, w)^{-1},
\label{eq:sumlog:gu}
\\
g_v &= -\zeta(u, v, w)^{-1} (\varphi \bigl( \tfrac{w}{v}) - d \bigr) - \tfrac{1}{v},
\label{eq:sumlog:gv}
\\
g_{w_i} &= -\zeta(u, v, w)^{-1} \tfrac{v}{w_i} - \tfrac{1}{w_i} & \forall i \in \iin{d}.
\label{eq:sumlog:gwi}
\end{align}
\end{subequations}

Note that $\bomega > 1$ since:
\begin{subequations}
\begin{align}
(p, q, r) &\in \intr(\klog^\ast)
\\
\Rightarrow p d - q + p \tsum_{i \in \iin{d}} \log \bigl(-\tfrac{r_i}{p} \bigr) & < 0
\\
\Leftrightarrow \tfrac{1}{d} \bigl(1 + d - \tfrac{q}{p} + \tsum_{i \in \iin{d}} \log \bigl(-\tfrac{r_i}{p} \bigr) \bigr) - \log(d) &> \tfrac{1}{d} + \log \bigl( \tfrac{1}{d} \bigr)
\label{eq:omegaboundc}
\\
\Leftrightarrow \omega \bigl( \tfrac{1}{d} \bigl(1 + d - \tfrac{q}{p} + \tsum_{i \in \iin{d}} \log \bigl( -\tfrac{r_i}{p} \bigr) \bigr) - \log(d) \bigr) &> \tfrac{1}{d}
\label{eq:omegaboundd}
\\
\Leftrightarrow \bomega &> 1,
\end{align}
\end{subequations}
where \cref{eq:omegaboundd} follows from \cref{eq:omegaboundc} by applying $\omega$ to both sides, and noting that $\omega$ is the inverse of $x \mapsto x + \log(x)$ on the reals (from the form of \cref{eq:omegachar}).

We would like to find $\tg \coloneqq (\tg_p, \tg_q, \tg_r)$ such that $-g(-\tg) = (p, q, r)$.
Fix $\zeta \coloneqq \zeta(-\tg_p, -\tg_q, -\tg_r)$.
Then, from \cref{eq:sumlog:gu}:
\begin{equation}
    p = -\zeta^{-1}.
    \label{eq:sumlog:p}
\end{equation}
Combining \cref{eq:sumlog:gwi} with \cref{eq:sumlog:p}, we need for all $i \in \iin{d}$:
\begin{equation}
    r_i = \tfrac{\tg_q}{\zeta \tg_{r_i}} - \tfrac{1}{\tg_{r_i}}
    = -\tfrac{p \tg_q}{\tg_{r_i}} - \tfrac{1}{\tg_{r_i}}
    \Rightarrow \tg_{r_i} = \tfrac{-p \tg_q - 1}{r_i}
    .
\label{eq:sumlog:tgr}
\end{equation}
Combining \cref{eq:sumlog:gv} with \cref{eq:sumlog:p,eq:sumlog:tgr}, we need:
\begin{equation}
    q = \zeta^{-1} \bigl( \varphi \bigl(\tfrac{\tg_r}{\tg_q} \bigr) - d \bigr) - \tfrac{1}{\tg_q}
    = -p \big( \varphi \big(\tfrac{-p \tg_q - 1}{r \tg_q} \big) - d \big) - \tfrac{1}{\tg_q}.
\label{eq:sumlogq}
\end{equation}
Replacing the definition of $\varphi$ in \cref{eq:sumlogq} and rearranging:
\begin{subequations}
\begin{align}
    \tfrac{q}{p} - d - \tsum_{i \in \iin{d}} \log(-\tfrac{r_i}{p}) - 1 &= -\tfrac{1}{p \tg_q} - d \log(1 + \tfrac{1}{p \tg_q}) - 1
    \\
    \Rightarrow -\tfrac{q}{p} + d + \tsum_{i \in \iin{d}} \log(-\tfrac{r_i}{p}) + 1 &= 1 + \tfrac{1}{p \tg_q} + d \log(1 + \tfrac{1}{p \tg_q}).
    \label{eq:omegac}
\end{align}
\end{subequations}
Note \cref{eq:omegac} has the form $\beta = a + d \log(a)$.
Letting $a = d b$:
\begin{equation*}
    \beta = d b + d \log(d b) 
    \Rightarrow \tfrac{\beta}{d} - \log(d) = b + \log(b).
\end{equation*}
Therefore,
\begin{subequations}
\begin{align}
    \tfrac{1}{d} \bigl( 1 + d - \tfrac{q}{p} + \tsum_{i \in \iin{d}} \log(-\tfrac{r_i}{p}) \bigr) - \log(d) &= 1 + \tfrac{1}{p \tg_q} + \log \bigl(1 + \tfrac{1}{p \tg_q} \bigr)
    \label{eq:omegad}
    \\
    \Rightarrow \bomega &= 1 + \tfrac{1}{p \tg_q}
    \\
    \Rightarrow \tg_q &= -p^{-1} (1 - \bomega)^{-1}.
    \label{eq:tgqfinal}
\end{align}
\end{subequations}
Substituting \cref{eq:tgqfinal} in \cref{eq:sumlog:tgr} gives \cref{eq:tgsumlog:ri}.
Finally, due to \cref{eq:grdot}:
\begin{equation*}
    \tg_p = \tfrac{-d - 2 - q \tg_q - \langle r, \tg_r \rangle}{p}
    = \tfrac{-d - 2 + \sfrac{q}{p} + 2 \bomega}{p (1 - \bomega)}.
\end{equation*}\qed
\subsection{Hypograph Power Cone}
\label{sec:conjproofs:hpower}

{\it Proof of \cref{lem:hpower}}
For $\varphi$ given in \cref{sec:hpower}, let $\hat{y}$ denote the root of $h$ from \cref{lem:hpower},
which we show is unique.
Note that $\hat{y}$ must satisfy:
\begin{equation*}
    \varphi(\hat{y} e - p \alpha_i) = \varphi(r)
    \Rightarrow \varphi ( \sfrac{\hat{y}}{\alpha} - p e ) = \varphi ( \sfrac{r}{\alpha} ).
\end{equation*}
Since $(p, r) \in \intr(\khpower^\ast)$, this implies $\hat{y} > 0$.
The derivatives of $h$ are:
\begin{align*}
    h'(y) &= \tsum_{i \in \iin{d}} \tfrac{\alpha_i}{y - p \alpha_i},
    \\
    h''(y) &= \tsum_{i \in \iin{d}} -\tfrac{\alpha_i}{(y - p \alpha_i)^2} < 0.
\end{align*}
Note that $p < 0$ for $(p, r) \in \intr(\khpower^\ast)$ and therefore $h'(y) > 0$ for all $y > \max_{i \in \iin{d}} \{ p \alpha_i \}$, i.e. the domain of $h$.
So the root of $h$ is unique.
Since $h$ is concave and increasing, a root-finding Newton-Raphson method will converge quadratically from any initial $y_{-} < \hat{y}$ (\cref{lem:lr}).
We may pick, for example, $y_{-} = 0$, which ensures $y_{-} < \hat{y}$ and $y_{-}$ is in the domain of $h$. 
\qed

{\it Proof of \cref{prop:hpower}}
For convenience, let $\zeta$ be the function $\zeta(u, w) \coloneqq \varphi(w) - u$, where $\varphi$ is from \cref{sec:hpower}.
Let $(u, w) \in \intr(\khpower)$.
Then the gradient of $f$ with respect to components $u$ and $w$ is:
\begin{subequations}
\begin{align}
    g_u &= \zeta(u, w)^{-1},
    \label{eq:hpower:gu}
    \\
    g_{w_i} &= -\zeta(u, w)^{-1} \tfrac{ \varphi(w) \alpha_i}{w_i} - \tfrac{1}{w_i} & \forall i \in \iin{d}.
    \label{eq:hpower:gwi}
\end{align}
\end{subequations}
We would like to find $\tg \coloneqq \tg(p, r)$ such that $-g(-\tg) = (p, r)$.
Using \cref{eq:hpower:gu}, we need:
\begin{equation}
    p = -\zeta(-\tg_p, -\tg_r)^{-1} = -(\varphi(-\tg_r) + \tg_p)^{-1}.
    \label{eq:powermeanp}
\end{equation}
Using \cref{eq:hpower:gwi,eq:powermeanp}, we need for all $i \in \iin{d}$:
\begin{equation}
    r_i = -\zeta(-\tg_p, -\tg_r)^{-1} \tfrac{\varphi(-\tg_{r}) \alpha_i}{\tg_{r_i}} - \tfrac{1}{\tg_{r_i}}
    = -\tfrac{1 - p \alpha_i \varphi(-\tg_{r})}{\tg_{r_i}}.
\label{eq:powermeanr}
\end{equation}
From \cref{eq:powermeanp,eq:powermeanr}:
\begin{subequations}
\begin{align}
    \tg_p &= -p^{-1} - \varphi(-\tg_r),
    \\
    \tg_{r_i} &= \tfrac{p \alpha_i \varphi(-\tg_{r}) - 1}{r_i} & \forall i \in \iin{d}.
\end{align}
\label{eq:powermeantg}
\end{subequations}
It remains to show how to evaluate $\varphi(-\tg_{r})$.
Applying $\varphi$ from \cref{sec:hpower} to both sides of \cref{eq:powermeanr}, after collecting for each $i$:
\begin{subequations}
\begin{align}
    \varphi(r) &= \tprod_{i \in \iin{d}} \bigl( -{\tg_{r_i}}^{-1} \bigl( -p \alpha_i \varphi(-\tg_{r}) + 1 \bigr) \bigr)^{\alpha_i}
    \\
    \color{blue}
    &= \varphi(-{\tg_{r}})^{-1} \tprod_{i \in \iin{d}} \bigl( -p \alpha_i \varphi(-\tg_{r}) + 1 \bigr)^{\alpha_i}
    \\
    \color{black}
     &= \tprod_{i \in \iin{d}} \varphi(-{\tg_{r}})^{-\alpha_i} \bigl( -p \alpha_i \varphi(-\tg_{r}) + 1 \bigr)^{\alpha_i}
    \\
     &= \tprod_{i \in \iin{d}} \big( -p \alpha_i + \varphi(-\tg_{r})^{-1} \big)^{\alpha_i}
    \\
    \Rightarrow \log(\varphi(r)) &= \tsum_{i \in \iin{d}} \alpha_i \log \bigl( -p \alpha_i + \varphi(-\tg_{r})^{-1} \bigr)
    \\
    \Rightarrow h(\varphi(-\tg_{r})^{-1}) &= 0.
    \label{eq:geomphi}
\end{align}
\end{subequations}
From \cref{eq:geomphi}, we can evaluate $\varphi(-\tg_{r})^{-1}$ easily due to \cref{lem:hpower}.
Combining this with \cref{eq:powermeantg} justifies \cref{eq:khpower:tg}.
In the special case where $\alpha = \sfrac{e}{d}$, we can solve $h(y) = 0$ exactly, giving, $\varphi(-\tg_{r})^{-1} = (\varphi(r) + \sfrac{p}{d})^{-1}$.
Substituting this in \cref{eq:powermeantg} gives \cref{eq:khgeom:tg}.
\qed
\subsection{Radial Power Cone}
\label{sec:conjproofs:rpower}

{\it Proof of \cref{lem:rpower}}
Let $\hat{y}$ denote the positive root of $h$ from \cref{lem:rpower}, which we show is unique.
In the special case where $\alpha = \sfrac{e}{d_2}$, we can solve $h(y) = 0$ exactly.
It can be verified that $\hat{y} > 0$ is given by:
\begin{equation}
    \hat{y} = -p^{-1} + d_2 \tfrac{p + \sqrt{\varphi(r) (\sfrac{d_2^2}{p^2} \varphi(r) + d_2^2 - 1)}}{\varphi(r) d_2^2 - p^2},
    \label{eq:rpowerequal}
\end{equation}
where $\varphi$ is given in \cref{sec:rpower}.
Note that the denominator is positive since $(p, r) \in \intr(\krgeom^\ast)$.
Let us turn to the case of non-uniform $\alpha$.
The first two derivatives of $h$ are:
\begin{align*}
    h'(y) &= 2 \tsum_{i \in \iin{d_2}} \tfrac{\alpha_i^2}{\alpha_i y + \sfrac{(1 + \alpha_i)}{p}} - 2 \tfrac{y + \sfrac{1}{p}}{y (y + \sfrac{2}{p})},
    \\
    h''(y) &= -2 \tsum_{i \in \iin{d_2}} \tfrac{\alpha_i^3}{( \alpha_i y + \sfrac{(1 + \alpha_i)}{p} )^2} + \tfrac{2 (y^2 + \sfrac{2}{p^2} )}{y^2 (y + \sfrac{2}{p} )^2}.
\end{align*}
We have that $h$ is decreasing for $p, y > 0$ (and cannot have more than one root), since:
\begin{equation}
    h'(y) < 2 \tsum_{i \in \iin{d_2}} \tfrac{\alpha_i^2}{a_i y + \sfrac{(1 + \alpha_i)}{p}} - 2 \tfrac{y}{y (y + \sfrac{2}{p} )}
    \leq 2 \tfrac{1}{y + \sfrac{(1 + 1)}{p}} - 2 \tfrac{1}{y + \sfrac{2}{p}} = 0.
\label{eq:rpower:hp}
\end{equation}
In the second inequality we use the fact that the preceding expression is convex in $\alpha$, and therefore maximized at an extreme point of the simplex that $\alpha$ belongs to.
Similar reasoning shows that $h''(y) > 0$ for $p, y > 0$:
\begin{equation}
    h''(y) > -2\tsum_{i \in \iin{d_2}} \tfrac{\alpha_i^3}{( \alpha_i y + \sfrac{(1 + \alpha_i)}{p})^2} + \tfrac{2}{(y + \sfrac{2}{p})^2}
    \geq -\tfrac{2}{(y + \sfrac{2}{p} )^2} + \tfrac{2}{ (y + \sfrac{2}{p})^2} = 0.
\label{eq:rpower:hpp}
\end{equation}
The second inequality follows from the fact that the preceding expression is concave in $\alpha$.

Due to \cref{eq:rpower:hp,eq:rpower:hpp}, a root-finding Newton-Raphson method converges quadratically starting from some $y_{-}$ such that $y_{-}  \leq \hat{y}$ (\cref{lem:lr}).
We may use the solution from the equal powers case, i.e. \cref{eq:rpowerequal} for $y_{-}$.
To see why, note that the function $\bar{h}(\alpha, y) = \tsum_{i \in \iin{d_2}} 2 \alpha_i \log(2 \alpha_i y^2 + (1 + \alpha_i) \sfrac{2 y}{p} ) - \log(\varphi(r)) - \log(\sfrac{2 y}{p} + y^2) - 2 \log(\sfrac{2 y}{p})$ is convex and symmetric in $\alpha$ (ignoring $y$, it can be checked that the Hessian of $\bar{h}$ is diagonal with nonnegative entries).
So for any fixed $y$, $\bar{h}(\alpha, y)$ is minimized at $\alpha = \sfrac{e}{d_2}$.
Since $h$ is decreasing, a solution to $\bar{h}(\sfrac{e}{d_2}, y) = 0$ lower bounds $\hat{y}$.
\qed

{\it Proof of \cref{prop:rpower}}
For convenience, let $\zeta$ be the function $\zeta(u, w) \coloneqq \varphi(w) - \lVert u \rVert^2$, where $\varphi$ is from \cref{sec:rpower}.
Let $(u, w) \in \intr(\krpower)$.
Then the gradient of $f$ with respect to components $u$ and $w$ is:
\begin{subequations}
\begin{align}
    g_{u_i} &= \tfrac{2 u_i}{\zeta(u, w)} & \forall i \in \iin{d_1},
    \label{eq:gu}
    \\
    g_{w_i} &= \tfrac{-2 \alpha_i \varphi(w)}{w_i \zeta(u, w)} - \tfrac{1 - \alpha_i}{w_i} & \forall i \in \iin{d_2}.
    \label{eq:gw}
\end{align}
\label{eq:rpowergrad}
\end{subequations}
We would like to find $\tg \coloneqq \tg(p, r)$ such that $-g(-\tg) = (p, r)$.
First, if $p = 0$, it is easy to see from \cref{eq:rpowergrad} that:
\begin{align*}
    \tg_{p_i} &= 0 & \forall i \in \iin{d_1},
    \\
    \tg_{r_i} &= -\tfrac{1 + \alpha_i}{r_i} & \forall i \in \iin{d_2}.
\end{align*}

For the case $p \neq 0$, let us show that without loss of generality, we may assume $p \in \bbR_+$ (and therefore $d_1 = 1$).
Fix $\zeta \coloneqq \zeta(-\tg_p, -\tg_r)$.
Let $Q \in \bbR^{d_1 \times d_1}$ be a suitable Householder transformation mapping $p \in \bbR^{d_1}$ to a vector of zeros except for one entry that is equal to $\lVert p \rVert$.
Since the function $\zeta$ is invariant to orthonormal transformations on the first input, $f(u, w) = f(Q u, w)$.
It is also easy to see from the definition of the dual gradient in \cref{eq:gpdef} that this implies:
\begin{equation}
    \tg_p = Q^\top \tg_{Q p}.
    \label{eq:Qgfact}
\end{equation}
From \cref{eq:gu}:
\begin{align}
    \tg_p = \tfrac{\zeta p}{2}.
    \label{eq:rpower:tgpmulti}
\end{align}
Let $\tg_{\lVert p \rVert}$ denote the $p$-component of the conjugate gradient at $(\lVert p \rVert, r) \in \intr(\krpower^\ast)$ with $d_1 = 1$.
Due to \cref{eq:Qgfact} and the invariance of $\zeta$ to transformation by $Q^\top$:
\begin{equation*}
    \zeta = \zeta(-\tg_p, -\tg_r) = \zeta(- Q^\top \tg_{Q p}, -\tg_r) = \zeta(-\tg_{\lVert p \rVert}, -\tg_r) = 
    \tfrac{2 \tg_{\lVert p \rVert}}{\lVert p \rVert}.
\end{equation*}
Substituting into \cref{eq:rpower:tgpmulti}:
\begin{equation*}
    \tg_p = \tfrac{\tg_{\lVert p \rVert} \cdot p}{\lVert p \rVert},
\end{equation*}
where $\tg_{\lVert p \rVert} \in \bbR_+$ can be computed as for the $p \in \bbR_+$ case.

Suppose $p \in \bbR_+$ from hereon.
Due to \cref{eq:gu}, we need:
\begin{align}
    \zeta &= \tfrac{2 \tg_p}{p}.
    \label{eq:rpower:zeta}
\end{align}
Since $\zeta > 0$, we have $\tg_p > 0$. From \cref{eq:gw}:
\begin{subequations}
\begin{align}
    r_i &= \tfrac{2 \alpha_i \varphi(-\tg_r)}{-\tg_{r_i} \zeta} - \tfrac{1 - \alpha_i}{\tg_{r_i}} & \forall i \in \iin{d_2},
    \label{eq:powria}
    \\
    \Rightarrow r_i (-\tg_{r_i}) \zeta &= 2 \alpha_i \varphi(-\tg_r) + (1 - \alpha_i) \zeta & \forall i \in \iin{d_2}.
    \label{eq:powri}
\end{align}
\end{subequations}
Applying $\varphi$ to both sides, across all $i$:
\begin{equation*}
    \varphi(r) \varphi(-\tg_{r}) \zeta^2 = \varphi(2 \alpha_i \varphi(-\tg_r) + (1 - \alpha_i) \zeta).
\end{equation*}
Substituting for $\varphi(-\tg_r) = \zeta + {\tg_p}^2 = \tfrac{2 \tg_p}{p} + {\tg_p}^2$ and the expression for $\zeta$ from \cref{eq:rpower:zeta}:
\begin{equation}
    \varphi(r) \bigl(\tfrac{2 \tg_p}{p} + {\tg_p}^2 \bigr) \bigl( \tfrac{2 \tg_p}{p} \bigr)^2 = \varphi \big( 2 \alpha_i (\tfrac{2 \tg_p}{p} + {\tg_p}^2) + (1 - \alpha_i) \tfrac{2 \tg_p}{p} \big).
    \label{eq:rpowerroot}
\end{equation}
We treat this as a root-finding problem for $\tg_p$.
Taking the logarithm of both sides in \cref{eq:rpowerroot} and rearranging, we would like to find a root for $h(\tg_p) = 0$, where $p, \tg_p > 0$.
This can be solved easily by \cref{lem:rpower}.
Due to \cref{eq:powri}, the solution can be used to compute $\tg_r$ with:
\begin{align}
    \tg_{r_i} &= -\tfrac{1}{r_i \zeta} ( 2 \alpha_i \varphi(-\tg_r) + (1 - \alpha_i) \zeta)
    = -\tfrac{\alpha_i (1 + p \tg_p) + 1}{r_i} &
    \forall i \in \iin{d_2}
    ,
    \label{eq:rpowerproof:tgri}
\end{align}
which shows \cref{eq:rpower:tg}.
Combining \cref{eq:rpowerequal} with \cref{eq:rpowerproof:tgri} gives \cref{eq:krgeom:tg}.
\qed
\subsection{Infinity Norm Cone}
\label{sec:conjproofs:hpower:inf}

{\it Proof of \cref{lem:inf}}
Let $\hat{y}$ denote the negative root of $h$ from \cref{lem:inf}, which we show is unique.
Note that if $r = 0$, then trivially $\hat{y} = -\sfrac{(d + 1)}{p}$ ($p \neq 0$ on the $\intr(\kinf^\ast)$).
So assume that $r \neq 0$.
The derivatives of $h$ are:
\begin{align*}
    h'(y) &= p + y \tsum_{i \in \iin{d}} r_i^2 (1 + r_i^2 y^2)^{-1/2},
    \\
    h''(y) 
    &= \tsum_{i \in \iin{d}} r_i^2(1 + r_i^2 y^2)^{-3/2} > 0.
\end{align*}
Once again, $h$ can have at most one root $\hat{y}$ on the halfline $y \leq 0$ since:
\begin{equation*}
    h'(y) \geq p + y \tsum_{i \in \iin{d}} \tfrac{r_i^2}{\sqrt{r_i^2 y^2}}
    = p - \tsum_{i \in \iin{d}} \lvert r_i \rvert
    > 0.
\end{equation*}
The last inequality follows from $(p, r) \in \intr(\kinf^\ast)$.
Since $h$ is increasing and convex, a root-finding Newton-Raphson method will converge quadratically from any $y_{+} \geq \hat{y}$ (\cref{lem:lr}).
Consider the function:
\begin{equation*}
    \bar{h}(y) \coloneqq y \bigl( p - \tsum_{i \in \iin{d}} \lvert r_i \rvert \bigr) + 1
    = p y + \tsum_{i \in \iin{d}} \lvert r_i \rvert \lvert y \rvert + 1
    \leq h(y).
\end{equation*}
The root of $\bar{h}$ is at $y = -(p - \tsum_{i \in \iin{d}} \lvert r_i \rvert)^{-1}$.
Since $h$ is increasing in $y$, we may use this root for $y_{+}$.
Alternatively, we could use $\bar{h}(y) \coloneqq p y + d + 1 \leq h(y)$, and its root gives $y_{+} = -\sfrac{(d + 1)}{p}$.
\qed

{\it Proof of \cref{prop:inf}}
Let $(u, w) \in \intr(\kinf)$ and define $\zeta_i(u, w) \coloneqq u^2 - w_i^2$ for all $i \in \iin{d}$.
Then the gradient of $f$ is:
\begin{subequations}
\begin{align}
    g_u &= \tfrac{d - 1}{u} - \tsum_{i \in \iin{d}} \tfrac{2 u}{\zeta_i(u, w)},
    \label{eq:infu}
    \\
    g_{w_i} &= \tfrac{2 w_i}{\zeta_i(u, w)} & \forall i \in \iin{d}.
    \label{eq:infw}
\end{align}
\end{subequations}
We would like to find $\tg \coloneqq \tg(p, r)$ such that $-g(-\tg(p, r)) = (p, r)$.
Let us fix $\zeta_i \coloneqq {\tg_p}^2 - {\tg_{r_i}}^2$ for all $i \in \iin{d}$.
From \cref{eq:infw}, for all $i \in \iin{d}$:
\begin{equation*}
    r_i = \tfrac{2 \tg_{r_i}}{\zeta_i}
    \Rightarrow \tfrac{1}{2} ({\tg_p}^2 - {\tg_{r_i}}^2) r_i = \tg_{r_i}.
\end{equation*}
This implies the signs of $r_i$ and $\tg_{r_i}$ equal for all $i \in \iin{d}$, and:
\begin{equation}
    \tg_{r_i} = 1_{r_i \neq 0} \cdot \tfrac{\sqrt{1 + r_i^2 {\tg_p}^2} - 1}{r_i}.
    \label{eq:epinorminf:tgri}
\end{equation}
Substituting into the definition for $\zeta_i$, for all $i \in \iin{d}$:
\begin{equation}
    \zeta_i = 
    1_{r_i = 0} \cdot {\tg_p}^2 + 1_{r_i \neq 0} \cdot \tfrac{-2 + 2 \sqrt{1 + r_i^2 {\tg_p}^2}}{r_i^2}.
    \label{eq:epinorminf:zi}
\end{equation}
From \cref{eq:infu,eq:epinorminf:zi}:
\begin{subequations}
\begin{align}
    d - 1 - \tsum_{i \in \iin{d}} \tfrac{2 {\tg_p}^2}{\zeta_i} &= p \tg_p
    \\
    \Leftrightarrow d - 1 - \tsum_{i \in \iin{d} : r_i \neq 0} \tfrac{2 {\tg_p}^2 r_i^2}{-2 + 2 \sqrt{1 + r_i^2 {\tg_p}^2}} - \tsum_{i \in \iin{d} : r_i = 0} 2 &= p \tg_p
    \\
    \Leftrightarrow d - 1 - \tsum_{i \in \iin{d} : r_i \neq 0} \tfrac{{\tg_p}^2 r_i^2 (-1 - \sqrt{1 + r_i^2 {\tg_p}^2})}{1 - (1 + r_i^2 {\tg_p}^2)} - \tsum_{i \in \iin{d} : r_i = 0} 2 &= p \tg_p
    \\
    \Leftrightarrow d - 1 + \tsum_{i \in \iin{d}} \tfrac{{\tg_p}^2 r_i^2 (-1 - \sqrt{1 + r_i^2 {\tg_p}^2})}{r_i^2 {\tg_p}^2} &= p \tg_p
    \\
    \Leftrightarrow p \tg_p + \tsum_{i \in \iin{d}} \sqrt{1 + r_i^2 {\tg_p}^2} + 1 &= 0.
    \label{eq:epinorminf:tgp}
\end{align}
\end{subequations}
We treat \cref{eq:epinorminf:tgp} as a root-finding problem in the variable $\tg_p < 0$, which can be easily solved due to \cref{lem:inf}.
Expression \cref{eq:inf:tgr} is obtained from \cref{eq:epinorminf:tgri}.
\qed
\section{Proofs of Inverse Hessian Operators}
\label{sec:hessproofs}

We reuse the notation from \cref{sec:invhess}.
For convenience, we let $\tp(t) = \tp + t \tx$ and we let $\tg(t)$ denote the conjugate gradient at $\tp(t)$.
We use $'$ to denote derivatives with respect to the linearization variable $t$, i.e. $H^\ast(\tp)\cdot \tx =  {\tg}'(0)$.
\subsection{Hypograph Power Cone}
\label{sec:hessproofs:hpower}

{\it Proof of \cref{prop:invhess:hpower}}
Let us differentiate \cref{eq:powermeanp} and \cref{eq:powermeanr} at $\tp(t)$ with respect to $t$.
Differentiating \cref{eq:powermeanp}:
\begin{subequations}
\begin{align}
    p + t x &= -(\varphi(-\tg_r(t)) + \tg_p(t))^{-1}
    \\
    \Rightarrow x &= (\varphi(-\tg_r(t)) + \tg_p(t))^{-2} (\tfrac{\diff}{\diff t} \varphi(-\tg_r(t)) + {\tg_p}'(t)).
    \label{eq:xi}
\end{align}
\end{subequations}
Differentiating \cref{eq:powermeanr}, for all $i \in \iin{d}$:
\begin{subequations}
\begin{align}
    r_i + t z_i &=
    -\tfrac{1 - p(t) \alpha_i \varphi(-\tg_{r}(t))}{\tg_{r_i}(t)}
    \\
    \begin{split}
    \Rightarrow z_i &= \tg_{r_i}(t)^{-2} {\tg_{r_i}}'(t) (1 - p(t) \alpha_i \varphi (-\tg_r(t)))
    + \tg_{r_i}(t)^{-1} \alpha_i (x  \varphi(-\tg_r(t)) + {}
    \\
    &\pheq p(t) \tfrac{\diff}{\diff t} \varphi(-\tg_r(t))).
    \end{split}
    \label{eq:zi}
\end{align}
\end{subequations}
Now \cref{eq:xi,eq:zi} give a nonlinear system that we wish to solve for ${\tg}'(t)$.
From the chain rule:
\begin{equation*}
    \tfrac{\diff}{\diff t} \varphi(-\tg_r(t)) 
    = \varphi(-\tg_r(t)) \tsum_{i \in \iin{d}} \alpha_i (\tg_{r_i}(t))^{-1} {\tg_{r_i}}'(t).
\end{equation*}
From hereon let us drop the variable $t$ from our notation for brevity.
Define:
\begin{equation*}
    K \coloneqq \tsum_{i \in \iin{d}} \alpha_i {\tg_{r_i}}^{-1} {\tg_{r_i}}'.
\end{equation*}
From \cref{eq:zi}, for all $i \in \iin{d}$:
\begin{equation}
    {\tg_{r_i}}' 
    = {\tg_{r_i}}^{2} \cdot \tfrac{z_i + {\tg_{r_i}}^{-1} \alpha_i (-x  \varphi(-\tg_r) - p \varphi(-\tg_r) K )}{1 -p \alpha_i \varphi (-\tg_r)}.
    \label{eq:grip}
\end{equation}
Multiplying \cref{eq:grip} by $\alpha_i {\tg_{r_i}}^{-1}$ and summing over all $i \in \iin{d}$ gives:
\begin{align*}
    K 
    &= \tsum_{i \in \iin{d}} \tfrac{z_i \alpha_i \tg_{r_i} + \alpha_i^2 \varphi(-\tg_r) (-x - p K )}{-p \alpha_i \varphi (-\tg_r) + 1}
    \\
    \Rightarrow K &= \big( 1 + \tsum_{j \in \iin{d}} \tfrac{ \alpha_j^2 \varphi(-\tg_r) p}{1 - p \alpha_j \varphi (-\tg_r)} \big)^{-1} \tsum_{i \in \iin{d}} \tfrac{z_i \alpha_i \tg_{r_i} + \alpha_i^2 \varphi(-\tg_r) (-x)}{1 - p \alpha_i \varphi (-\tg_r)}.
\end{align*}
Let $k_{1,i} \coloneqq 1 - \alpha_i \varphi(-\tg_r) p$ for all $i \in \iin{d}$, $k_2 \coloneqq \tsum_{i \in \iin{d}} \sfrac{\alpha_i^2}{k_{1,i}}$, and $k_3 \coloneqq 1 + \varphi(-\tg_r) p k_2$.
Then from \cref{eq:xi,eq:grip}:
\begin{subequations}
\begin{align}
    {\tg_p}' &= (\varphi(-\tg_r) + \tg_p)^2 x - K \varphi(-\tg_r)
    \\
    &= \bigl((\varphi(-\tg_r) + \tg_p)^2 + \tfrac{k_2 {\tg_p}^2}{k_3} \bigr) x - \tfrac{\varphi(-\tg_r)}{k_3} \langle z, \tfrac{\alpha \tg_{r}}{k_1}  \rangle,
    \\
    {\tg_{r_i}}' &= {\tg_{r_i}}^{2} \cdot \tfrac{z_i + {\tg_{r_i}}^{-1} \alpha_i (-x  \varphi(-\tg_r) - p \varphi(-\tg_r) K )}{k_{1,i}}
    \\
    &= \tfrac{{\tg_{r_i}}^2}{k_{1,i}} z_i - \tfrac{\alpha_i \tg_{r_i}}{k_{1,i}} \tfrac{\varphi(-\tg_r)}{k_3} x - \tfrac{p \varphi(-\tg_r)}{k_3} \langle z, \tfrac{\alpha \tg_{r}}{k_1}  \rangle \tfrac{\alpha_i \tg_{r_i}}{k_{1,i}} \quad i \in \iin{d}.
\end{align}
\label{eq:hasthpower}
\end{subequations}
Now the desired result \cref{eq:hinvhpower} can be obtained by invoking \cref{eq:dirHfact}.
Recall that in \cref{eq:hasthpower}, ${\tg}' = H^\ast(\tp) \cdot \tx$.
So we simply \emph{replace} instances of $\tp$ and $\tg$ in \cref{eq:hasthpower} by $-g(\tu)$ and $-\tu$ respectively to obtain \cref{eq:hinvhpower}.
\qed
\subsection{Radial Power Cone}
\label{sec:hessproofs:rpower}

{\it Proof of \cref{prop:invhess:rpower}}
Let:
\begin{equation}
    K \coloneqq \tsum_{i \in \iin{d_2}} \alpha_i (\tg_{r_i}(t))^{-1} {\tg_{r_i}}'(t).
    \label{eq:invhess:rpower:K}
\end{equation}
For $\varphi$ defined in \cref{sec:rpower}, we have that:
\begin{align*}
    \tfrac{\diff}{\diff t} \varphi(-\tg_r(t)) &= \tsum_{i \in \iin{d_2}} \varphi(-\tg_r(t)) \cdot 2 \alpha_i (-\tg_{r_i}(t))^{-1} (-{\tg_{r_i}}'(t))
    \\
    &= 2 \varphi(-\tg_r) K.
\end{align*}

We fix $\zeta \coloneqq \varphi(-\tg_r(t)) - \lVert \tg_p(t) \rVert^2$ for convenience.
We start by rearranging and differentiating \cref{eq:rpower:zeta} at $\tp(t)$ with respect to a linearization variable $t$, which we omit for brevity:
\begin{subequations}
\begin{align}
    x_i &= \tfrac{\diff}{\diff t} \bigl( \tfrac{2 \tg_{p_i}(t)}{\zeta(\tp(t))} \bigr)
    \\
    &= 2 \big( \tfrac{{\tg_{p_i}}'}{\zeta} - \tfrac{\tg_{p_i}}{\zeta^2} (2 \varphi(-\tg_r)K - 2 \langle \tg_p, {\tg_p}' \rangle) \big) & \forall i \in \iin{d_1}.
    \label{eq:rpowerxj}
\end{align}
\end{subequations}

We use \cref{eq:rpowerxj} to obtain a new expression for $\langle \tg_p, {\tg_p}' \rangle$:
\begin{subequations}
\begin{align}
    \langle x, \tg_p \rangle &= 2 \Bigl( \tfrac{\langle {\tg_p}', \tg_p \rangle}{\zeta} - \tfrac{\lVert \tg_p \rVert^2}{\zeta^2} (2 \varphi(-\tg_r)K - 2 \langle \tg_p, {\tg_p}' \rangle) \Bigr)
    \\
    \Rightarrow \langle \tg_p, {\tg_p}' \rangle &= \tfrac{\zeta^2}{2 k_1}  \bigl( \langle x, \tg_p \rangle + \tfrac{4 \lVert \tg_p \rVert^2}{\zeta^2} \varphi(-\tg_r)K \bigr)
    ,
    \label{eq:tgpinphi2}
\end{align}
\end{subequations}
where $k_1 \coloneqq \varphi(-\tg_r) + \lVert \tg_p \rVert^2$.
Next, we differentiate \cref{eq:powri}, for all $i \in \iin{d_2}$:
\begin{equation*}
\begin{split}
    -z_i \tg_{r_i} \zeta - r_i {\tg_{r_i}}' \zeta - 2 r_i \tg_{r_i} ( \varphi(-\tg_r) K - \langle \tg_p, {\tg_p}' \rangle)
    = {}
    \\
    \quad 4 \alpha_i \varphi(-\tg_r) K + 2 (1 - \alpha_i) ( \varphi(-\tg_r) K - \langle \tg_p, {\tg_p}' \rangle),
\end{split}
\end{equation*}
whence, for all $i \in \iin{d_2}$ (replacing $\langle \tg_p, {\tg_p}' \rangle$ using \cref{eq:tgpinphi2}):
\begin{subequations}
\begin{align}
    -r_i {\tg_{r_i}}' \zeta
    &= 4 \alpha_i \varphi(-\tg_r) K + 2 (1 - \alpha_i + r_i \tg_{r_i}) ( \varphi(-\tg_r) K - \langle \tg_p, {\tg_p}' \rangle) + z_i \tg_{r_i} \zeta 
    \\
    &= 4 \alpha_i \varphi(-\tg_r) K - 4 \alpha_i \tfrac{\varphi(-\tg_r)}{\zeta} ( \varphi(-\tg_r) K - \langle \tg_p, {\tg_p}' \rangle) + z_i \tg_{r_i} \zeta
    \\
    &=  \tfrac{2 \alpha_i \varphi(-\tg_r)}{k_1} \big( 2 \lVert \tg_p \rVert^2 K + \zeta \langle \tg_p, x \rangle \big) + z_i \tg_{r_i} \zeta
    .
    \label{eq:powgrip}
\end{align}
\end{subequations}
In the second equality we use the identity, for all $i \in \iin{d_2}$:
\begin{equation*}
    1 - \alpha_i + r_i \tg_{r_i} = -2 \alpha_i \tfrac{\varphi(-\tg_r)}{\zeta},
\end{equation*}
which can be derived as follows.
From \cref{eq:powria}, we have that for all $i \in \iin{d_2}$:
\begin{subequations}
\begin{align}
    r_i \tg_{r_i} &= -\tfrac{2 \alpha_i \varphi(-\tg_r)}{\zeta} - (1 - \alpha_i)
    \label{eq:rpowertricka}
    \\
    \Rightarrow \tsum_{j \in \iin{d_2} : j \neq i} -r_j \tg_{r_j} \alpha_i &= 
    \tsum_{j \in \iin{d_2} : j \neq i} \Bigl( \tfrac{2 \alpha_i \alpha_j \varphi(-\tg_r)}{\zeta} + 
    \alpha_i (1 - \alpha_j) \Bigr).
    \label{eq:rpowertrickb}
\end{align}
\end{subequations}
Multiplying \cref{eq:rpowertricka} by $(1 - \alpha_i)$ and adding \cref{eq:rpowertrickb}, we get:
\begin{align}
    (1 - \alpha_i) r_i \tg_{r_i} - \tsum_{j \in \iin{d_2} : j \neq i} \alpha_i r_j \tg_{r_j} &= -1 + d_2 \alpha_i & \forall i \in \iin{d_2}.
    \label{eq:powertrickc}
\end{align}
Combining \cref{eq:powertrickc} with \cref{eq:grdot}, for all $i \in \iin{d_2}$:
\begin{subequations}
\begin{align}
    r_i \tg_{r_i} + \tfrac{-1 + d_2 \alpha_i - (1 - \alpha_i) r_i \tg_{r_i}}{-\alpha_i} + \langle p, \tg_p \rangle &= -d_2 - 1
    \\
    \Rightarrow 1 - \alpha_i + r_i \tg_{r_i} &= \alpha_i \bigl( \langle p, -\tg_p \rangle - 2 \bigr) .
    \label{eq:rigrifact}
\end{align}
\end{subequations}
Combining \cref{eq:rigrifact} with \cref{eq:rpower:zeta} gives the result:
\begin{equation*}
    1 - \alpha_i + r_i \tg_{r_i} = -2 \alpha_i \Bigl( \tfrac{\lVert \tg_p \rVert^2}{\zeta} + 1 \Bigr)
    = -2 \alpha_i \tfrac{\varphi(-\tg_r)}{\zeta}.
\end{equation*}

Next, substituting into ${\tg_{r_i}}'$ from \cref{eq:powgrip} into \cref{eq:invhess:rpower:K}:
\begin{align*}
    K &= -\tsum_{i \in \iin{d_2}} \alpha_i (\tg_{r_i})^{-1} \tfrac{2 \alpha_i \varphi(-\tg_r) k_1^{-1} \big( 2 \lVert \tg_p \rVert^2 K + \zeta \langle \tg_p, x \rangle \big) - z_i (-\tg_{r_i}) \zeta }{r_i \zeta}
    \\
    &= \tfrac{-2 k_2 k_1^{-1} \varphi(-\tg_r) \langle x, \tg_p \rangle - \langle \sfrac{\alpha}{r}, z \rangle }{1 + 4 k_2 k_1^{-1} \varphi(-\tg_r) \lVert \tg_p \rVert^2},
\end{align*}
where $k_2 \coloneqq \langle \sfrac{\alpha}{r}, \sfrac{\alpha}{\tg_r} \rangle$.
To obtain ${\tg_p}'$, we replace $\langle \tg_p, {\tg_p}' \rangle$ in \cref{eq:rpowerxj} using the expression in  \cref{eq:tgpinphi2}, then substitute for $K$.
Define for convenience $k_3 \coloneqq \sfrac{k_1}{2 \varphi(-\tg_r)} + \sfrac{2 k_2 \lVert \tg_p \rVert^2}{\zeta}$.
Then:
\begin{subequations}
\begin{align}
    {\tg_p}' &= \zeta \big(\tfrac{x}{2} + \tfrac{\tg_p}{\zeta^2} (2 \varphi(-\tg_r)K - 2 \langle \tg_p, {\tg_p}' \rangle) \big)
    \\
    &= \tfrac{\zeta}{2} x - \tfrac{\tg_p}{k_3} \big( \tfrac{2 k_2 \varphi(-\tg_r) + \zeta k_3}{k_1} \langle x, \tg_p \rangle + \langle \sfrac{\alpha}{r}, z \rangle \big).
\end{align}
\label{eq:rpower:tgpp}
\end{subequations}

Finally, substituting for $K$ in \cref{eq:powgrip} and rearranging, for all $i \in \iin{d_2}$:
\begin{equation}
    {\tg_{r_i}}' = -\tfrac{\tg_{r_i}}{r_i} z_i + \tfrac{\alpha_i}{k_3 r_i} \bigr( - \langle x, \tg_p \rangle + \tfrac{2 \lVert \tg_p \rVert^2}{\zeta} \langle \sfrac{\alpha}{r}, z \rangle \bigr).
    \label{eq:rpower:tgrp}
\end{equation}
As in \cref{sec:invhess:hpower}, applying \cref{eq:dirHfact} to \cref{eq:rpower:tgpp,eq:rpower:tgrp} gives the result \cref{eq:Hinvrpower}.
\qed
%



\begin{acknowledgements}
This work has been partially funded by the National Science Foundation under grant OAC-1835443 and the Office of Naval Research under grant N00014-18-1-2079.
The authors thank the anonymous reviewers and editor for their comments and suggestions.
\end{acknowledgements}

\bibliography{main}

\end{document}

%% file: defs.tex
\DeclareMathOperator{\diff}{d}
\DeclareMathOperator{\cl}{cl}
\DeclareMathOperator{\intr}{int}
\DeclareMathOperator{\Diag}{Diag}

\DeclareMathOperator{\logdet}{logdet}
\DeclareMathOperator{\rtdet}{rtdet}
\DeclareMathOperator{\hpower}{hpower}
\DeclareMathOperator{\rpower}{rpower}

\DeclareMathOperator{\spec}{spec}
\DeclareMathOperator{\hgeom}{hgeom}
\DeclareMathOperator{\rgeom}{rgeom}

\DeclareMathOperator{\argsup}{argsup}

\newcommand{\bomega}{\bar{\omega}}
\newcommand{\bOmega}{\bar{\Omega}}

\newcommand{\tg}{g^\ast}
\newcommand{\tx}{\tilde{x}}

\newcommand{\tu}{\tilde{u}}
\newcommand{\tp}{\tilde{p}}
\newcommand{\tprod}{\textstyle\prod}
\newcommand{\tsum}{\textstyle\sum}
\newcommand{\K}{\mathcal{K}}
\newcommand{\bbR}{\mathbb{R}}
\newcommand{\bbS}{\mathbb{S}}

\newcommand{\iin}[1]{\llbracket #1 \rrbracket}

\newcommand{\klog}{\K_{\log}}
\newcommand{\klogdet}{\K_{\logdet}}
\newcommand{\khgeom}{\K_{\hgeom}}
\newcommand{\krgeom}{\K_{\rgeom}}
\newcommand{\khpower}{\K_{\hpower}}
\newcommand{\krpower}{\K_{\rpower}}
\newcommand{\krtdet}{\K_{\rtdet}}
\newcommand{\kinf}{\K_{\ell_{\infty}}}
\newcommand{\kspec}{\K_{\ell_{\spec}}}

\newcommand{\pheq}{\hphantom{{}={}}}

%% file: abstract.tex
The recent interior point algorithm by Dahl and Andersen \cite{dahl2021primal} for nonsymmetric cones as well as earlier works \cite{nesterov2012towards,nesterov1999infeasible} require derivative information from the conjugate of the barrier function of the cones in the problem.
Besides a few special cases, there is no indication of when this information is efficient to evaluate.
We show how to compute the gradient of the conjugate barrier function for seven useful nonsymmetric cones.
In some cases this is helpful for deriving closed-form expressions for the inverse Hessian operator for the primal barrier.

%% file: introduction.tex
\section{Introduction}
\label{sec:introduction}

Many convex optimization problems can be expressed as the minimization of a linear function over an intersection of \emph{symmetric} cones.
These include the nonnegative, (rotated) second order, and positive semidefinite cones.
Many popular solvers for programming with symmetric cones implement primal-dual interior point methods (PDIPMs) specialized for these cones only (e.g. \cite{borchers1999csdp,yamashita2003implementation,vandenberghe2010cvxopt}).
These specialized PDIPMs enjoy properties such as efficiently computable \emph{symmetric} search directions (that ensure the algorithm is invariant to exchanging of roles of the primal and dual programs), which usually rely on forming so-called \emph{scaling matrices} for symmetric cones, first suggested by \cite{nesterov1997self}.
In addition, the symmetric cones are self-dual and admit barrier functions with highly efficient oracles for PDIPMs.
Unfortunately, many properties of symmetric cones are not straightforward to generalize for other conic sets.
In spite of the fact that many useful sets can only be modeled directly or can be modeled more efficiently with cones that are not symmetric, there are significantly fewer implementations of algorithms for \emph{nonsymmetric} cones.

Early algorithmic frameworks for nonsymmetric cones have been suggested by \cite{nesterov1999infeasible} and \cite{nesterov2012towards}, but have not been implemented by any solver we are aware of.
The high-level algorithms in both papers require information relating to the conjugates of the barrier functions of the cones in the problem, and it is not clear when this can be evaluated efficiently.
Specifically, the only nonsymmetric cones with known, efficient procedures for evaluating the conjugate of the barrier function (or its gradient) are the three-dimensional exponential cone \cite{serrano2015algorithms}, the three-dimensional power cone \cite{nesterov2012towards}, and the cone of sparse positive semidefinite matrices with chordal sparsity \cite{andersen2013logarithmic}.
Additionally, \cite{nesterov1999infeasible} requires solving a linear system in each iteration that is twice the size of the linear systems that arise in symmetric algorithms.

More recent nonsymmetric algorithms bypass some of the aforementioned issues.
These include the algorithm by Skajaa and Ye \cite{skajaa2015homogeneous}, which has been implemented in the solvers Alfonso \cite{papp2021alfonso} and Hypatia \cite{coey2020solving}.
A key advantage of this algorithm is that it requires very few oracles for each primal cone in the problem, and doesn't require oracles relating to the conjugate of the primal barrier.

In this paper we are primarily motivated by the alternative algorithm of Dahl and Andersen \cite{dahl2021primal}, which is implemented by the MOSEK solver \cite{mosek2020modeling}.
The algorithm is based on a technique by \cite{tunccel2001generalization,myklebust2014interior} which generalizes the concept of scaling matrices from symmetric cones.
Like the algorithm of \cite{skajaa2015homogeneous}, the linear systems solved in each iteration are equal in size to those arising in symmetric algorithms.
Unlike the algorithm of \cite{skajaa2015homogeneous}, the gradient of the conjugate barrier is required to compute search directions.
However, the search directions satisfy a number of desirable properties such as ensuring that the violations of residuals and complementarity conditions decrease at the same rate.
The authors also propose a neighborhood that allows stepping further away from the central path compared to \cite{skajaa2015homogeneous}, although there is no proof of polynomial time convergence.
One might expect, owing to the use of conjugate barrier information, that the search directions of \cite{dahl2021primal} would allow convergence in fewer iterations than the search directions of \cite{skajaa2015homogeneous}.

The search directions in the algorithm by \cite{dahl2021primal} relate to those from symmetric algorithms as follows.
Given a primal-dual pair of points $(w, r)$ for a symmetric cone, there exists a unique scaling matrix $T$ satisfying the secant equations $r = T w$ and $g(w) = T g^\ast(r)$, where $g$ and $g^\ast$ denote the gradients of the primal barrier and its conjugate \cite{nesterov1997self,nesterov1998primal}.
A key idea of \cite{myklebust2014interior,dahl2021primal} is to construct a general positive definite matrix for any cone, satisfying the two secant equations.
The authors choose a specific formula for this scaling matrix that requires calculating the Hessian of a primal barrier function, and adding BFGS updates to the Hessian.
These low rank updates include adding multiples of the conjugate gradient.

The conjugate gradient can always be evaluated via a numerical procedure (e.g. applying Newton's method to an optimization problem), but this generic approach is computationally slow and can be numerically challenging.
In particular, it requires applying the inverse Hessian of the primal barrier in each iteration (which can become a new bottleneck in a PDIPM).
A large number of damped Newton iterations may be necessary to get near the region of quadratic convergence towards the end of the PDIPM, when the distance to the cone boundary is small.
Our aim is to show efficient methods of calculating conjugate gradients for seven useful nonsymmetric cones, which an interior point solver could support.
Aside from their use in the algorithm by \cite{dahl2021primal}, procedures for evaluating conjugate gradients are useful due to their applications in the frameworks by \cite{nesterov1999infeasible,nesterov2012towards}.
The cones for which we are able to demonstrate efficient procedures for evaluating conjugate gradients are:
\begin{description}[itemsep=0pt]
\item[{\bf Logarithm cone,}]
the hypograph of the perspective of the sum of logarithms.
\item[{\bf Log-determinant cone,}]
the hypograph of the perspective of the log-deter-minant function.
\item[{\bf Hypograph power cone,}]
the hypograph of the power mean function.
\item[{\bf Root-determinant cone,}]
the hypograph of the $d$th-root-determinant function.
\item[{\bf Radial power cone,}]
the generalized power cone by \cite[Section 3.1.2]{chares2009cones}.
\item[{\bf Infinity norm cone,}]
the epigraph of the $\ell_\infty$ norm function.
\item[{\bf Spectral norm cone,}]
the epigraph of the spectral norm function.
\end{description}
In all cases, computing conjugate gradients can be reduced to simple numerical procedures (e.g. univariate root-finding). 
For some special cases (e.g.  where the hypograph or radial power cones are parametrized by equal powers), the conjugate gradients can be written as closed-form expressions. 

Our conjugate gradient evaluations lead to a second result.
By differentiating the procedure to calculate the conjugate gradients we are able to derive closed-form expressions for the inverse Hessian of the primal barrier.
This is useful for measuring central path proximity and some of the linear system methods suggested by \cite{coey2021performance}.
The inverse Hessian is already known for five of the cones above \cite{coey2021conic,coey2022thesis}, so we only show it for the hypograph and radial power cones.

%% file: preliminaries.tex
\section{Preliminaries}

\subsection{Notation}

We use $\bbR_+^d$, $\bbR_{++}^d$, $\bbR_{-}^d$, and $\bbR_{--}^d$ to denote the $d$-dimensional vectors of nonnegative, positive, nonpositive, and negative reals respectively.
Likewise, $\bbS_+^d$ and $\bbS_{++}^d$ are the positive semidefinite and positive definite matrices respectively with side dimension $d$.
For a natural number $d$, we define the index set $\iin{d} \coloneqq \{1, 2, \ldots, d\}$.
$1_{P}$ is the indicator function, that is equal to one if statement $P$ is true and zero otherwise.
For a set $C$, $\cl(C)$ and $\intr(C)$ denote the closure and interior of $C$ respectively.
We use round parentheses, e.g. $(a, b, c)$ for Cartesian products of vectors or matrices, $\langle \cdot, \cdot \rangle$ for an inner product between vectors or matrices, and $e$ to denote a vector of ones.
$V^\top$ is the transpose of $V$ with respect to $\langle \cdot, \cdot \rangle$.
Throughout this paper, division and multiplication between vectors should be interpreted as componentwise.
The operator $\Diag$ maps a vector to a matrix with the vector on the diagonal.
We write $g(w)$ or $H(w)$ to denote gradients or Hessians evaluated at $w$, and use subscripts to refer to partial derivatives.

\subsection{Cones and Barrier Functions}

A proper cone $\K$ is closed, convex, pointed, and full-dimensional.
The \emph{dual cone} of $\K$ is $\K^\ast \coloneqq \{ r \in \bbR^d : \langle w, r \rangle \geq 0, \forall w \in \K \}$.
$\K^\ast$ is a proper cone if and only if $\K$ is a proper cone.

Analysis of conic interior point methods relies on associating with each cone a \emph{logarithmically-homogeneous, self-concordant barrier} (LHSCB).
Following \cite[Sections 2.3.1 and 2.3.3]{nesterov1994interior}, $f : \intr(\K) \to \bbR$ is a $\nu$-LHSCB for $\K$, where $\nu \geq 1$ is the \emph{LHSCB parameter}, if it satisfies $f(w_i) \to \infty$ along every sequence $w_i \in \intr(\K)$ converging to the boundary of $\K$, and:
\begin{subequations}
\begin{alignat}{3}
\big\lvert \nabla^3 f(w) [r, r, r] \big\rvert 
& \leq 2 \bigl( \nabla^2 f(w) [r, r] \bigr)^{3/2} 
&\quad&\forall w \in \intr(\K), r \in \bbR^d,
\label{eq:lhscb:sc}
\\
f(\theta w) 
& = f(w) - \nu \log(\theta) 
&\quad&\forall w \in \intr(\K), \theta \in \bbR.
\label{eq:lhscb:lh}
\end{alignat}
\end{subequations}
As a consequence of \eqref{eq:lhscb:lh}, the gradient $g$ of $f$ satisfies \cite[Equation (2.5)]{nesterov1997self}:
\begin{equation}
    \langle -g(w), w \rangle = \nu \quad \forall w \in \K.
    \label{eq:grdot}
\end{equation}
For the LHSCB $f$, we define the convex conjugate, $f^\ast$, as the function:
\begin{equation}
    f^\ast(r) \coloneqq \textstyle\sup_{w \in \intr(\K)} \{ - \langle r, w \rangle - f(w) \}.
    \label{eq:fast}
\end{equation}
In fact, $f^\ast$ is an LHSCB for $\K^\ast$ \cite[Equation (2.6)]{nesterov1997self} and we refer to it as the \emph{conjugate barrier}.
The gradient $g^\ast$ of $f^\ast$ may be defined from the unique solutions to the optimization problem in \cref{eq:fast}:
\begin{equation}
    \tg(r) = -\argsup_{w \in \intr(\K)} \{ - \langle r, w \rangle - f(w) \}.
    \label{eq:gpdef}
\end{equation}
From \cref{eq:fast}, \cref{eq:gpdef} and \cref{eq:grdot}:
\begin{equation}
    f^\ast(r) = -\langle r, -\tg(r) \rangle - f(-\tg(r)) = -\nu - f(-\tg(r)).
    \label{eq:fast2}
\end{equation}
We refer to $g^\ast$ as the \emph{conjugate gradient} (which has no relation to the method of conjugate gradients).
The negative gradients of LHSCBs $f$ and $f^\ast$ are bijective linear maps between $\K$ and $\K^\ast$.
In particular \cite[Theorem 2.5]{myklebust2014interior}:
\begin{equation}
    -g^\ast(-g(w)) = w \quad \forall w \in \K, \quad
    -g(-g^\ast(r)) = r \quad \forall r \in \K^\ast.
\label{eq:grbilinear}
\end{equation}
Note that \cref{eq:grbilinear} characterizes the gradient and conjugate gradient maps as negative inverses of each other.

Some of the LHSCBs we use are related to \emph{unitarily invariant} functions \cite{lewis1995convex}.
Let $W = U_W \Diag(\sigma_W) V_W^\top$ be the singular value decomposition of $W \in \bbR^{d_1 \times d_2}$, such that $\sigma_W \in \bbR^{d_1}$.
If $W$ is symmetric, then $\sigma_W$ are the eigenvalues of $W$ and $U_W = V_W$.
Suppose $F : \bbR^{d_1 \times d_2} \to \bbR$ is a function given by $F(W) = f(\sigma_W)$, where $f \in \bbR^{d_1} \to \bbR$ is some symmetric function (invariant to the order of its inputs).
Then $F$ is unitarily invariant.
Let $G$ and $g$ denote the gradients of $F$ and $f$.
In \cref{sec:sumlog,sec:hpower,sec:inf}, we use \cite[Theorem 3.1]{lewis1995convex}:
\begin{equation}
    G(W) = U_W \Diag(g(\sigma_W)) V_W^\top.
    \label{eq:unitarily:G}
\end{equation}